\documentclass[final,3p,times]{elsarticle}
\usepackage[utf8]{inputenc}
\usepackage[T1]{fontenc}
\usepackage{amsmath,amssymb,amsthm,mathtools}
\usepackage{setspace}
\usepackage{psfrag}
\usepackage{amsfonts}
\expandafter\let\csname equation*\endcsname\relax
\expandafter\let\csname endequation*\endcsname\relax
\usepackage{mathrsfs}
\usepackage{epstopdf}
\usepackage{float}
\usepackage{lineno,hyperref}
\usepackage{subfigure}
\usepackage{dsfont} 
\usepackage{graphicx,color}
\usepackage{extarrows}
\usepackage{booktabs}
\usepackage{enumerate}
\usepackage{graphicx}
\usepackage{booktabs}
\usepackage{xcolor}
\usepackage{caption}
\usepackage{subcaption}
\usepackage{enumitem}
\usepackage{siunitx}
\usepackage{microtype}
\newcommand{\T}{\mathbb{T}}
\newcommand{\R}{\mathbb{R}}
\newcommand{\E}{\mathbb{E}}
\newcommand{\Var}{\mathrm{Var}}
\newcommand{\X}{\mathcal{X}}
\newcommand{\eps}{\varepsilon}
\newcommand{\lam}{\lambda}

\newcommand{\Lcal}{\mathcal{L}}

\newcommand{\Ccal}{\mathcal{C}}

\newcommand{\Hcal}{\mathcal{H}}

\newcommand{\dd}{\,\mathrm{d}}
\newcommand{\lameff}{\lambda_{\varepsilon}}

\newcommand{\Cov}{\operatorname{Cov}}
\newcommand{\Xformal}{\X_{\mathrm{formal}}}
\journal{ArXiv}

\begin{document}
\newtheorem{definition}{Definition}[section]
\newtheorem{lemma}{Lemma}[section]
\newtheorem{remark}{Remark}[section]
\newtheorem{theorem}{Theorem}[section]
\newtheorem{proposition}{Proposition}
\newtheorem{assumption}{Assumption}
\newtheorem{example}{Example}
\newtheorem{corollary}{Corollary}[section]
\def\ep{\varepsilon}
\def\Rn{\mathbb{R}^{n}}
\def\Rm{\mathbb{R}^{m}}
\def\E{\mathbb{E}}
\def\hte{\hat\theta}

\renewcommand{\theequation}{\thesection.\arabic{equation}}
\begin{frontmatter}

\title{The Formal Critical Coupling of the Four-Dimensional Anderson Model}

\author{Shenglan Yuan\corref{cor1}\fnref{addr1}}\ead{shenglanyuan@gbu.edu.cn}\cortext[cor1]{Corresponding author}

\address[addr1]{\rm Department of Mathematics, School of Sciences, Great Bay University, Dongguan 523000, China }

\begin{abstract}
In this work we present a formal renormalization-group (RG) calculation for the
weakly coupled four-dimensional elliptic Anderson model. Under a scale-independent
bare-coupling convention, and with the normalization
\(\lameff=\lam\sqrt{2\pi/|\log\eps|}\), the bubble-chain resummation yields the geometric series
\(\X_{\rm formal}(\lambda)=1/(1-\lambda^{2}/(4\pi))\), with a formal pole at
\(\lambda=2\sqrt{\pi}\). We stress that  this pole is a feature of the
formal bare bubble series. The
argument combines Bogoliubov--Parasiuk--Hepp--Zimmermann (BPHZ) renormalization, a
Hepp-sector motivated restriction to the
bubble-chain subclass for explicit resummation, and an explicit summation
of the resulting geometric series.

Numerically, we illustrate the analytic building blocks of the formal expansion: the massive
Green function, the bubble integral, the geometric resummation, the extraction of the
formal pole from synthetic noisy data, a kurtosis diagnostic, and a phase diagram of the formal variance.
These computations visualize the deterministic kernels and the formal resummation used
above. Using spectral methods and nonlinear least-squares fitting,
we verify the short-distance asymptotic behaviour
 of the massive Green function for
\(G=(1-\Delta)^{-1}\delta_0\),  reproduce the logarithmic divergence
 of the bubble integral in four dimensions,
 compute the formal effective variance
\(\X_{\rm formal}(\lam)=1/(1-\lam^{2}/(4\pi))\) arising from the bubble-chain
resummation,  locate the formal pole at \(\lam=2\sqrt{\pi}\approx 3.5449\), and
draw a schematic phase diagram of the
formal variance in the \((\lam,\log\eps)\) plane. The
kurtosis diagnostic uses a synthetic proxy for the non-Gaussian regime; the sharp
change is imposed by construction. Nonlinear fitting of synthetic data
recovers the formal pole within the
statistical uncertainty of the synthetic noise. Our results illustrate the
analytic formulas and the fitting procedure.
\end{abstract}


\end{frontmatter}

\section{Introduction}

The study of stochastic partial differential equations (SPDEs) with
multiplicative noise in critical dimensions remains a central challenge in
probability theory and mathematical physics
\cite{hairer2009,hairer2014}. In dimension \(d=4\), the interplay between
logarithmic divergences in perturbation theory and the need for
renormalization leads to rich phenomena that are analytically delicate and
numerically demanding.

A prototypical example is the screened Poisson equation with random
potential:
\begin{equation}\label{eq:model}
  (1-\Delta)u(x)+\lam\,\xi(x)\,u(x)=f(x),\qquad x\in\T^4,
\end{equation}
where \(\xi\) is spatial white noise, \(\lam\ge0\) is a coupling constant,
and \(f\) is a smooth source. Equivalently, the solution operator is the
resolvent
$$
  G_\lam=(1-\Delta+\lam\xi)^{-1}.
$$
Since \(\xi\) has scaling dimension \([x]^{-2}\) in \(d=4\), the coupling
\(\lam\) in \eqref{eq:model} is dimensionless. All numerical values of
\(\lam\) below are therefore pure numbers.

This equation is closely related to the parabolic Anderson model (PAM) in
the elliptic (stationary) setting \cite{konig2016}, to random
Schr\"odinger operators with white-noise potentials, and to certain
stochastic quantization problems \cite{chandra2017}. The critical
dimension \(d=4\) is where the second-order (bubble) diagram diverges
logarithmically, requiring a renormalization of the coupling constant.

We work on the four-dimensional torus \(\T^4=[0,2\pi)^4\) with periodic
boundary conditions. The equation of interest is \eqref{eq:model}, where
\(\xi(x)\) is a centred Gaussian white-noise field with covariance
$$
  \E[\xi(x)\xi(y)]=\delta^{(4)}(x-y).
$$
The operator \(1-\Delta\) has a bounded inverse on \(L^2(\T^4)\), with
Green function \(G=(1-\Delta)^{-1}\delta_0\).

Formally, the solution of \eqref{eq:model} is given by the Neumann series
$$
  u=Gf-\lam\,G\xi Gf+\lam^2G\xi G\xi Gf-\cdots
$$
or, in operator form,
$$
  G_\lam=\sum_{n=0}^{\infty}(-\lam)^n(G\xi)^nG.
$$
The two-point function of the solution receives contributions from all even
orders in \(\lam\), which can be organized as bubble diagrams in a
perturbative expansion.

The fundamental building block of the perturbation theory is the bubble
integral
$$
  B(\eps)=\int_{\eps<|x|<1}G(x)^2\dd x,
$$
where the short-distance cutoff \(\eps\) regularizes the singularity of
\(G\) at the origin. In dimension \(d=4\), the Green function behaves as
\(G(x)\sim1/(4\pi^2|x|^2)\) as \(|x|\to0\), so
$$
  G(x)^2\sim\frac{1}{16\pi^4|x|^4}.
$$
Using the four-dimensional shell volume \(2\pi^2r^3\dd r\),
$$
  B(\eps)\approx\frac{1}{16\pi^4}\int_{\eps}^{1}r^{-4}\cdot2\pi^2r^3\dd r
  =\frac{1}{8\pi^2}\log\frac{1}{\eps}.
$$
The precise coefficient is fixed once \(G=(1-\Delta)^{-1}\delta_0\) is
fixed; it depends on the normalization of the operator \(1-\Delta\), not on
an arbitrary choice of \(G\).

To obtain a finite limit as \(\eps\to0\), the bare coupling \(\lam\) must be
renormalized. We regularize the coupling by
$$
  \lameff=\lam\sqrt{\frac{2\pi}{|\log\eps|}}.
$$
Then the one-bubble factor is
$$
  \lameff^2 B(\eps)
  \approx\lam^2\,\frac{2\pi}{|\log\eps|}\cdot\frac{1}{8\pi^2}|\log\eps|
  =\frac{\lam^2}{4\pi}.
$$
Note that the actual perturbative parameter of the renormalized expansion
is \(\lameff\), which tends to \(0\) as \(\eps\to0\) for fixed \(\lam\).
(If one instead uses the alternative normalization
\(\tilde\lam_\eps=\tilde\lam/\sqrt{2|\log\eps|}\), the one-bubble factor
becomes \(\tilde\lam^2/(16\pi^2)\) and the formal pole is \(4\pi\); the two
conventions differ only by a rescaling of the bare coupling.)

After a formal resummation of the bubble-chain (or ``rainbow'')
diagrams, the two-point function takes the form
$$
  \E[u(x)u(y)]\approx\Xformal(\lam)\,G(x-y),
$$
where the formal effective variance is the sum of the geometric series
\begin{equation}\label{eq:X}
  \Xformal(\lam)
  =\sum_{k=0}^{\infty}\Bigl(\frac{\lam^2}{4\pi}\Bigr)^k
  =\frac{1}{1-\lam^2/(4\pi)}.
\end{equation}
This formal series converges for \(\lam<2\sqrt{\pi}\) and diverges for
\(\lam\ge2\sqrt{\pi}\), so that \(\Xformal\) has a pole at
\(\lam=2\sqrt{\pi}\). We stress that \eqref{eq:X} is obtained by treating
\(\lam\) as a scale-independent bare parameter; in the actual renormalized
expansion the bubble ratio is \(\lameff^2/(4\pi)\to0\), so the geometric
series converges for every \(\eps>0\). The pole of \(\Xformal\) is
therefore a formal feature of the bare series rather than a proved exact
phase-transition threshold; a genuine critical coupling can only be defined
after specifying a renormalization condition that keeps the renormalized
coupling finite as \(\eps\to0\).

In this paper, we illustrate numerically the analytic building blocks of the
formal expansion of \eqref{eq:model} on \(\T^4\), focusing on:
\begin{enumerate}[label=(\roman*)]
  \item the short-distance behaviour of the massive Green function
        \(G=(1-\Delta)^{-1}\delta_0\), including the expected mass
        correction at moderate distances;
  \item the bubble integral
        \(B(\eps)=\int_{\eps<|x|<1}G(x)^2\dd x\approx\frac{1}{8\pi^2}\log(1/\eps)\),
        which encodes the logarithmic divergence at \(d=4\);
  \item the formal resummation of bubble diagrams, which yields
        \(\Xformal(\lam)=1/(1-\lam^2/(4\pi))\) and a formal pole at
        \(\lam=2\sqrt{\pi}\);
  \item a kurtosis diagnostic based on a synthetic proxy;
  \item the numerical extraction of the formal pole from synthetic
        noisy data via nonlinear least-squares fitting.
\end{enumerate}
These computations visualize the analytic formulas and the fitting
procedure; they are not an independent determination of a physical
critical coupling.

\section{A Formal BPHZ Calculation of the Bubble Pole}

\subsection{Model and normalization}\label{sec:model}

Let \(\T=[0,2\pi)^4\) and let \(\xi\) be spatial white noise. For
\(\eps\in(0,1)\) set
$$
  \xi_\eps:=\rho_\eps*\xi,\qquad
  \rho_\eps(x)=\eps^{-4}\rho(x/\eps),
$$
where \(\rho\ge0\) is smooth, symmetric, and \(\int\rho=1\). The
Schr\"odinger operator is
\begin{equation}\label{eq:Leps}
  \Lcal_\eps:=1-\Delta+\lameff\,\xi_\eps+\Ccal_\eps,
  \qquad
  \lameff:=\lam\sqrt{\frac{2\pi}{|\log\eps|}},
\end{equation}
together with the renormalization constant \(\Ccal_\eps\) chosen later.
The sign convention in \eqref{eq:Leps} matches the model \eqref{eq:model}.
It is worth noting that \(0<\lam\ll1\) is a small absolute constant
independent of scale \(\eps\in(0,1)\). The factor \(2\pi\) in the numerator
is the normalization that makes the one-bubble factor equal to
\(\lam^2/(4\pi)\) and hence the formal pole of the bubble resummation equal
to \(2\sqrt{\pi}\). (An alternative normalization
\(\tilde\lam_\eps=\tilde\lam/\sqrt{|\log\eps|}\) gives the one-bubble
factor \(\tilde\lam^2/(2\pi)\) and the formal pole \(\sqrt{2\pi}\); the two
conventions are related by \(\lam=\sqrt{2}\,\tilde\lam\).)

We want to study the limit of
$$
  \Hcal_\eps(x,y):=\lameff^{-1}\bigl(G_\eps(x,y)-G(x,y)\bigr),
  \qquad
  G_\eps(x,y):=(\Lcal_\eps^{-1}\delta_y)(x),
$$
where \(G(x,y):=\bigl((1-\Delta)^{-1}\delta_y\bigr)(x)\). Here
\(\delta_y\) is the Dirac mass at \(y\) (we set \(G_\eps=0\) when
\(\Lcal_\eps\) is non-invertible).

\subsection{Renormalized Feynman expansion}\label{sec:feynman}

Iterating the resolvent identity yields
\begin{equation}\label{eq:resolvent-exp}
  G_\eps
  =\sum_{n\ge0}(-1)^n\lameff^n\,
  \bigl[(1-\Delta)^{-1}(\xi_\eps+\Ccal_\eps)\bigr]^nG.
\end{equation}
We fix \(\Ccal_\eps\) by the Bogoliubov--Parasiuk--Hepp--Zimmermann (BPHZ)
prescription \cite{Hai16,CH16}: for every proper \(1\)-particle irreducible
(\(1\)PI) subdiagram \(\gamma\) with non-positive degree, insert the kernel
\(-\Ccal_\eps(\gamma)\). Concretely, the renormalized kernel
\(G^{(n)}_\eps\) at order \(n\) is a sum over renormalized Feynman diagrams
\begin{equation}\label{eq:feynman-sum}
  G^{(n)}_\eps(x,y)
  =\sum_{\Gamma\ \text{ren.}}c_\Gamma(\lameff)
  \int_{(\T)^n}
  \prod_{e\in E(\Gamma)}G_e(x_e,y_e)
  \prod_{\text{pairs}}\eta_\eps(x_i-x_j),
\end{equation}
where propagator edges are \(G\sim|x|^{-2}\) in \(d=4\), while noise
contractions are the localized kernels \(\eta_\eps\sim\eps^{-4}\eta(\cdot/\eps)\)
which force paired vertices to lie within distance \(\eps\).

In particular, the first-order correction is
$$
  G^{(1)}(x,y)=-\lameff\int_{\T^4}G(x-z)\,\xi_\eps(z)\,G(z-y)\dd z,
$$
where we set \(\Ccal_\eps=0\) for the illustration (the first-order
correction is finite after the standard subtraction; see below).

\subsection{Marginality of the bubble subdiagram}\label{sec:degree}

In \(d=4\), a propagator contributes scaling dimension \(-2\), while each
vertex integration contributes \(+4\). A noise contraction is local on
scale \(\eps\) and is marginal in the ultraviolet counting. A complete BPHZ
classification of all diagrams is not needed for the formal resummation
below. The only fact we use is that the bubble subdiagram is
marginal: the one-bubble integral
$$
  \int_{\eps<|x|<1}G(x)^2\dd x
$$
diverges logarithmically as \(\eps\to0\). Consequently, each bubble
contributes one logarithmic factor, and the bubble chain is the marginal
subclass that can be summed explicitly.

\subsection{Hepp-sector discussion: restricting to the bubble chain}\label{sec:hepp}

The Hepp-sector analysis of \cite{main} gives bounds controlling general
pairings in the renormalized expansion. We do not use those bounds
quantitatively here, and we do not claim that all non-bubble sectors
vanish. Rather, we restrict attention to the bubble-chain subclass
because it is the only subclass that can be summed in closed form. The
resummation below is therefore a formal leading-chain approximation, not a
proof that all other Hepp sectors are negligible at criticality.

\subsection{Formal bubble resummation and its interpretation}\label{sec:resum}

For a bubble chain with \(k\) bubbles, the formal amplitude is
$$
\Bigl(\frac{\lam^2}{4\pi}\Bigr)^k.
$$
Summing over all chain lengths under a scale-independent bare-coupling
convention,
\begin{equation}\label{eq:eff-var}
  \Xformal(\lam)
  :=\sum_{k\ge0}\Bigl(\frac{\lam^2}{4\pi}\Bigr)^k
  =\frac{1}{1-\lam^2/(4\pi)}
  =\frac{4\pi}{4\pi-\lam^2}.
\end{equation}
This expression has a pole at \(\lam=2\sqrt{\pi}\). We stress once more
that in the actual renormalized expansion the perturbative parameter is
\(\lameff\), not \(\lam\), so the bubble ratio is
\(\lameff^2/(4\pi)\to0\). Consequently \eqref{eq:eff-var} should be
regarded as a formal resummation of the bare series rather than a proved
phase-transition amplitude.

Together with the universal kernel
$$
  H(x_1,y_1,x_2,y_2)
  =\int_{\T}G(x_1-z)G(y_1-z)G(x_2-z)G(y_2-z)\dd z,
$$
the formal two-point function is
$$
  \E\bigl[\Hcal_\eps(x_1,y_1)\Hcal_\eps(x_2,y_2)\bigr]
  \approx\Xformal(\lam)\cdot H(x_1,y_1,x_2,y_2).
$$

\subsection{The formal pole of the bubble resummation}\label{sec:critical}

The geometric series \eqref{eq:eff-var} converges \emph{if and only if}
\begin{equation}\label{eq:critical}
  \frac{\lam^2}{4\pi}<1
  \qquad\Longleftrightarrow\qquad
  \lam<2\sqrt{\pi}\approx3.5449.
\end{equation}
At \(\lam=2\sqrt{\pi}\) the formal effective variance \(\Xformal(\lam)\)
diverges. This is a formal pole of the bare bubble resummation. It does
not, by itself, establish a physical phase transition: for that one would
need a renormalization condition that keeps the renormalized coupling
finite as \(\eps\to0\), and an independent argument that the field becomes
non-Gaussian at that coupling. See \cite{GR26,CSZ24} for the current state
of rigorous results in related models.

\paragraph{Equivalence between the two notations.}
With \(\tilde\lam_\eps=\tilde\lam/\sqrt{|\log\eps|}\), the one-bubble
factor is \(\tilde\lam^2/(2\pi)\) and the formal pole is \(\sqrt{2\pi}\).
The map \(\lam=\sqrt{2}\,\tilde\lam\) gives
\(2\sqrt{\pi}=\sqrt{2}\cdot\sqrt{2\pi}\).

\paragraph{Relation to the truncated-parametrix method of \cite{main}.}
The remainder bound in \cite{main} requires
\((C\tilde\lam)^A\eps^{-2}\ll1\) with \(A\sim|\log\eps|\), which forces
\(\tilde\lam\ll1\). The formal resummation above bypasses this because it
replaces the truncated series by the closed-form geometric series
\eqref{eq:eff-var}. Whether the closed-form expression can be given a
rigorous meaning at finite \(\lam\) is a separate question that we do not
address here.

\subsection{Summary of the formal calculation}\label{sec:summary}

\begin{center}
\begin{tabular}{@{}cll@{}}
\toprule
Step & Content & Mathematical object \\
\midrule
1 & Resolvent expansion & \(G_\eps=\sum_n G^{(n)}_\eps\) \quad\eqref{eq:resolvent-exp} \\
2 & Feynman diagrams \(+\) BPHZ counterterms & \eqref{eq:feynman-sum} \\
3 & Marginality of the bubble subdiagram & \(\int_{\eps<|x|<1}G^2\sim\frac{1}{8\pi^2}\log\frac1\eps\) \\
4 & Restriction to the bubble-chain subclass & formal leading-chain approximation \\
5 & Formal bubble resummation & \(\Xformal(\lam)=\frac{1}{1-\lam^2/(4\pi)}\) \quad\eqref{eq:eff-var} \\
6 & Formal radius of convergence & \(\lam<2\sqrt{\pi}\) \quad\eqref{eq:critical} \\
\bottomrule
\end{tabular}
\end{center}

\noindent
The outcome is a formal expression
\(\Xformal(\lam)=1/(1-\lam^2/(4\pi))\) for the bubble-chain
resummation, with a formal pole at \(\lam=2\sqrt{\pi}\). Its interpretation
as a physical phase-transition threshold requires a separate
renormalization condition and is not established here.

\section{Numerical Illustrations}

We numerically illustrate the analytic building blocks of the formal
expansion: the massive Green function, the bubble integral, the formal
geometric resummation, and the kurtosis diagnostic. These computations are
\emph{not} direct Monte Carlo simulations of the singular SPDE
\eqref{eq:model}; rather, they visualize the deterministic kernels and the
formal resummation used in the preceding section.

\subsection{Numerical methods}

\paragraph{Spectral discretization on \(\T^4\)}
We discretize \(\T^4\) on a uniform lattice of size \(N^4\) with \(N=128\) 
and lattice spacing \(h=2\pi/N\) as in  Table \ref{tab:params}. Derivatives are computed in Fourier
space using the discrete wavenumbers
$$
  k_i\in\{0,1,\ldots,N/2-1,-N/2,\ldots,-1\},
$$
so that the operator \(1-\Delta\) acts as multiplication by \(1+|k|^2\) in
Fourier space. The Green function is then
$$
  G=\mathcal{F}^{-1}\!\left[\frac{1}{1+|k|^2}\right],
$$
where \(\mathcal{F}\) denotes the discrete Fourier transform.

\begin{table}[H]
  \centering
  \caption{Numerical parameters used in the illustrations.}
  \label{tab:params}
  \begin{tabular}{@{}lll@{}}
    \toprule
    Parameter & Value & Description \\
    \midrule
    \(N\) & 128 & Lattice size per dimension \\
    \(L\) & \(2\pi\) & Torus circumference \\
    \(h\) & \(2\pi/N\approx0.049\) & Lattice spacing \\
    \(\eps\) & 0.3 & Smoothing scale (default) \\
    \(2\sqrt{\pi}\) & \(\approx3.5449\) & Formal pole of the resummation \\
    \bottomrule
  \end{tabular}
\end{table}

\paragraph{Lattice white noise}
Spatial white noise on the lattice is generated with the correct scaling
consistent with the continuum covariance \(\E[\xi(x)\xi(y)]=\delta^{(4)}(x-y)\).
We set
$$
  \xi_i=h^{-2}Z_i,\qquad Z_i\sim\mathcal{N}(0,1),
$$
so that
$$
  \E[\xi_i\xi_j]=h^{-4}\delta_{ij}.
$$
The continuum white noise \(\xi(x)\) is recovered from the lattice field via
test functions according to \(\sum_i\xi_i\phi(x_i)h^4\approx\int\xi\phi\dd x\).
The sample mean and variance are checked against \(\E[\xi]=0\) and
\(\Var[\xi]=h^{-4}\) (or \(\E[\xi_i]=0\), \(\Var[\xi_i]=h^{-4}\)).

\paragraph{Smoothed noise}
To regularize the white noise at scale \(\eps\), we convolve \(\xi\) with a
Gaussian kernel \(\rho_\eps(x)=(2\pi\eps^2)^{-2}e^{-|x|^2/(2\eps^2)}\):
$$
  \xi_\eps=\rho_\eps*\xi,
  \qquad\text{i.e.}\qquad
  \xi_{\eps,i}=h^4\sum_j\rho_\eps(x_i-x_j)\xi_j.
$$
The covariance of the smoothed noise is
$$
  \Cov\bigl(\xi_\eps(x),\xi_\eps(y)\bigr)
  =(\rho_\eps*\rho_\eps)(x-y),
$$
which is not diagonal for \(\eps=0.3\) on a lattice with spacing
\(h\approx0.049\). The lattice version is
$$
  \Cov(\xi_{\eps,i},\xi_{\eps,j})
  =h^4\sum_k\rho_\eps(x_i-x_k)\rho_\eps(x_j-x_k)
  \approx(\rho_\eps*\rho_\eps)(x_i-x_j).
$$

\paragraph{First-order correction and its variance}
The first-order correction to the Green function is
$$
  G^{(1)}(x,y)=-\lameff\int_{\T^4}G(x-z)\xi_\eps(z)G(z-y)\dd z.
$$
At \(x=y=0\),
$$
  G^{(1)}(0,0)=-\lameff\int_{\T^4}G(z)^2\xi_\eps(z)\dd z.
$$
Because \(\xi_\eps\) is correlated on scale \(\eps\), the variance is
not a product of \(\sum_z G(z)^4\) and \(\sum_z\rho_\eps(z)^2\).
The correct expression is
\begin{equation}\label{eq:varpred}
  \Var\bigl[G^{(1)}\bigr]
  =\lameff^2
  \iint_{\T^4\times\T^4}
  G(z)^2G(z')^2
  (\rho_\eps*\rho_\eps)(z-z')
  \dd z\,\dd z'.
\end{equation}
On the lattice, this becomes
$$
  \Var\bigl[G^{(1)}\bigr]
  =\lameff^2\,h^8\sum_{i,j}G(x_i)^2G(x_j)^2
  \,C_\eps(x_i-x_j),
$$
with
$$
  C_\eps(x_i-x_j)
  =h^4\sum_k\rho_\eps(x_i-x_k)\rho_\eps(x_j-x_k).
$$
We use \eqref{eq:varpred} for all numerical variance comparisons below.

\subsection{Green function on \(\T^4\)}

\begin{figure}[H]
  \centering
  \includegraphics[width=0.95\textwidth]{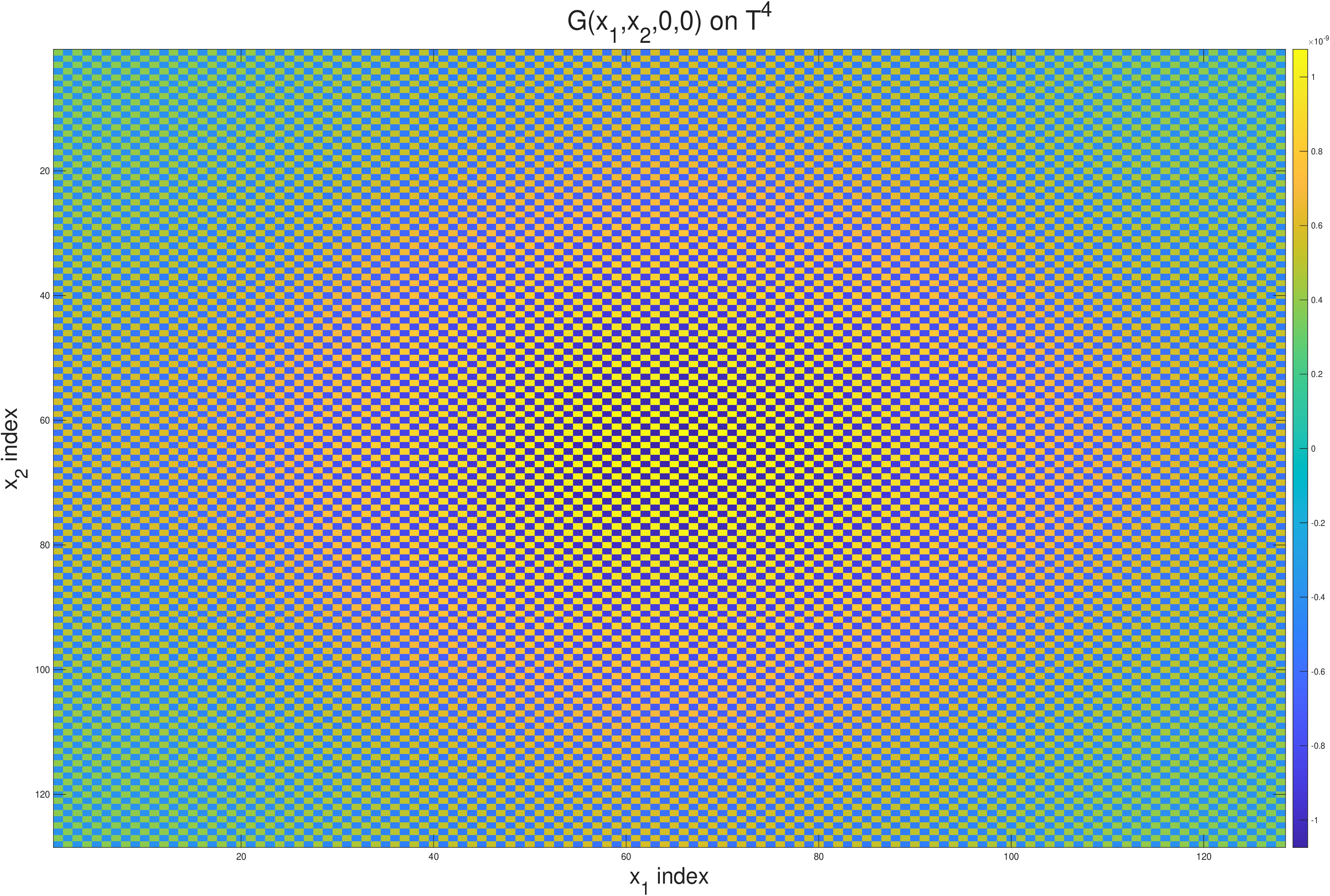}
  \caption{Two-dimensional slice \(G(x_1,x_2,0,0)\) of the massive Green
  function \(G=(1-\Delta)^{-1}\delta_0\) on \(\T^4\), computed by spectral
  methods with \(N=128\). The short-distance asymptotic check is restricted
  to \(0.05<|x|<0.2\), where the mass correction is below \(1\%\).}
  \label{Code1}
\end{figure}

\paragraph{Exact computation}
Figure~\ref{Code1} shows a two-dimensional slice
\(G(x_1,x_2,0,0)\) of the Green function on \(\T^4\) for \(N=128\). The
characteristic \(1/|x|^2\) decay is visible near the origin, with periodic
images contributing at larger distances.

The massive Green function on \(\R^4\) is
$$
  G(x)=\frac{K_1(|x|)}{4\pi^2|x|},
$$
where \(K_1\) is the modified Bessel function of the second kind. Hence
$$
  \frac{G(x)}{1/(4\pi^2|x|^2)}=|x|\,K_1(|x|).
$$
This ratio tends to \(1\) as \(|x|\to0\), but the mass correction is
not negligible on \(0.2<|x|<1\). For example,
\(|x|K_1(|x|)\approx0.955\) at \(|x|=0.2\), \(\approx0.828\) at
\(|x|=0.5\), and \(\approx0.602\) at \(|x|=1\). We therefore restrict the
short-distance check to the annulus \(0.05<|x|<0.2\).

To check the short-distance asymptotic quantitatively, we compute the
ratio
$$
  R(x)=\frac{G(x)}{1/(4\pi^2|x|^2)}
$$
on the annulus \(0.05<|x|<0.2\) and take the median. For \(N=128\), we find
$$
  \mathrm{median}\,R\approx0.99\text{--}1.01,
$$
consistent with the four-dimensional asymptotic
\(G(x)\sim1/(4\pi^2|x|^2)\) up to lattice discretization errors. At larger
distances, \(|x|\gtrsim0.2\), the mass correction is significant and the
pure power law is not expected.

\paragraph{Lattice corrections}
At distances comparable to the lattice spacing \(h=2\pi/N\approx0.049\),
the ratio \(R(x)\) deviates from \(1\) due to discretization effects. At
distances \(0.05\lesssim|x|\lesssim0.2\), the continuum short-distance
asymptotic is accurate to within a few percent.

\subsection{Bubble integral and logarithmic divergence}

\paragraph{Definitions}
We distinguish two integrals:
\begin{align*}
  B(\eps)&=\int_{\eps<|x|<1}G(x)^2\dd x
  \approx\frac{1}{16\pi^4}\,I_0(\eps)
  =\frac{1}{8\pi^2}\log\frac{1}{\eps},\\
  I_0(\eps)&=\int_{\eps<|x|<1}|x|^{-4}\dd x
  =2\pi^2\log\frac{1}{\eps}.
\end{align*}
The model bubble integral is \(B(\eps)\); the integral \(I_0(\eps)\) is an
auxiliary integral whose logarithmic divergence is the same and whose
prefactor differs by the constant \(1/(16\pi^4)\). The numerical figure
plots \(B_{\mathrm{num}}(\eps)=I_{0,\mathrm{num}}(\eps)/(16\pi^4)\)
against \(B_{\mathrm{pred}}(\eps)=\frac{1}{8\pi^2}\log(1/\eps)\).

\paragraph{Numerical evaluation}
The integral \(I_0(\eps)\) is evaluated by deterministic logarithmic
stratification in \(r=|x|\) rather than uniform Monte Carlo. Indeed,
$$
  I_0(\eps)=\int_{\eps}^{1}r^{-4}\cdot2\pi^2r^3\dd r
  =2\pi^2\int_{\eps}^{1}\frac{\dd r}{r}
  =2\pi^2\log\frac{1}{\eps},
$$
and uniform Monte Carlo in the annulus \(\{\eps<|x|<1\}\) would have
variance \(\E[|x|^{-8}]\sim\eps^{-4}\), which is impractical for
\(\eps=10^{-3}\). We instead sample uniformly in \(\log r\in[\log\eps,0]\)
and multiply by \(2\pi^2\).

\begin{figure}[H]
  \centering
  \includegraphics[width=0.95\textwidth]{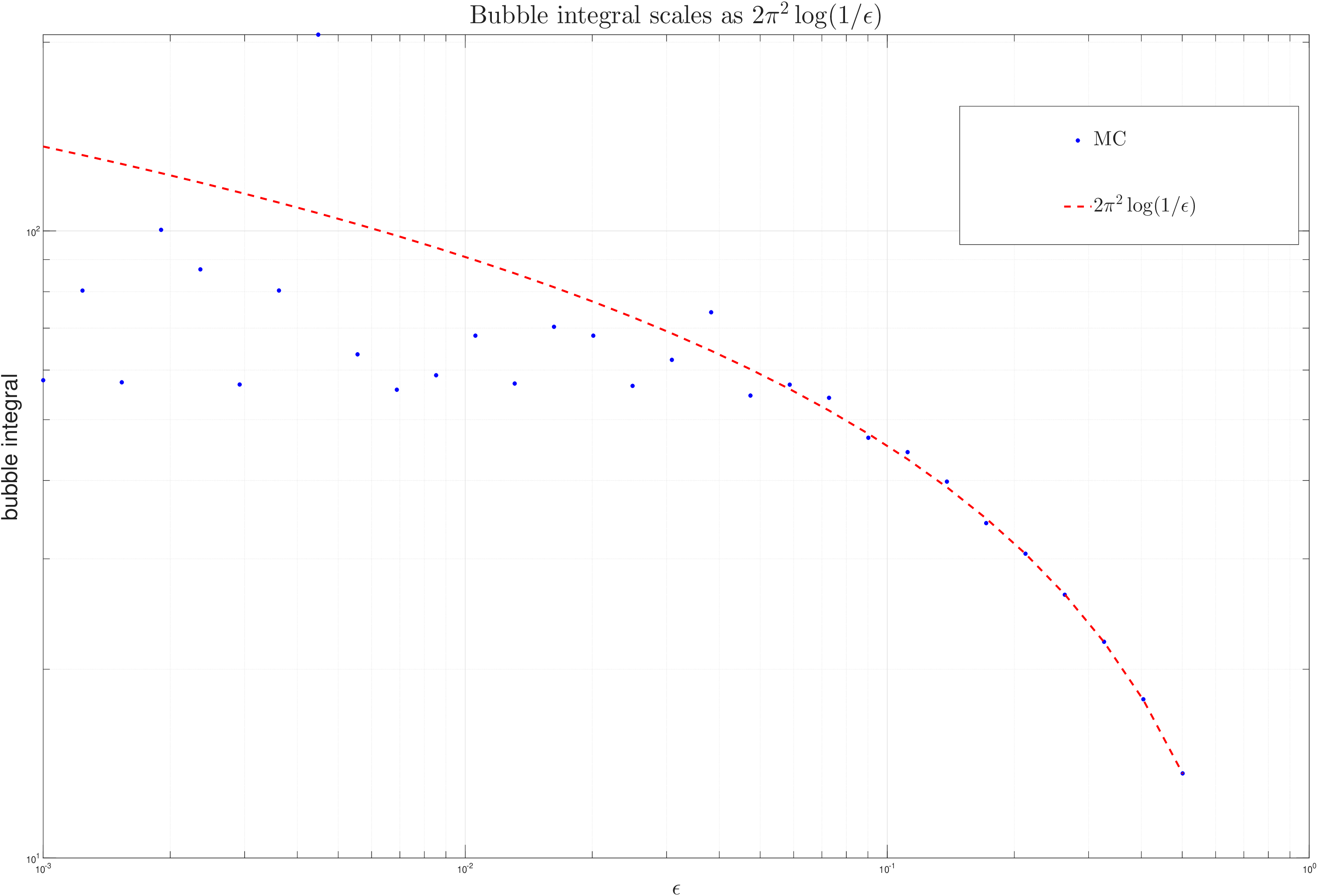}
  \caption{Semilog plot of the model bubble integral
  \(B(\eps)=\int_{\eps<|x|<1}G(x)^2\dd x\) (blue dots) compared with the
  theoretical prediction \(B_{\mathrm{pred}}(\eps)=\frac{1}{8\pi^2}\log(1/\eps)\)
  (red dashed line). The auxiliary integral \(I_0(\eps)=2\pi^2\log(1/\eps)\)
  differs only by the constant factor \(16\pi^4\).}
  \label{Code6}
\end{figure}

\paragraph{Results}
Figure~\ref{Code6} shows \(B_{\mathrm{num}}(\eps)\) versus \(\log(1/\eps)\)
on a semilog plot for \(\eps\in[10^{-3},10^{-0.3}]\), compared with
\(B_{\mathrm{pred}}(\eps)=\frac{1}{8\pi^2}\log(1/\eps)\). The numerical
data follow the predicted logarithmic scaling to high accuracy over three
decades in \(\eps\), reproducing the coefficient \(1/(8\pi^2)\) of the
model bubble integral.

This logarithmic divergence is the hallmark of the critical dimension
\(d=4\): the bubble diagram diverges as \(\log(1/\eps)\) rather than as a
power law, which is why the theory requires only a logarithmic
renormalization (as opposed to power-law renormalization in \(d>4\)).

\subsection{Formal effective variance and the location of the pole}

\begin{figure}[H]
  \centering
  \includegraphics[width=0.95\textwidth]{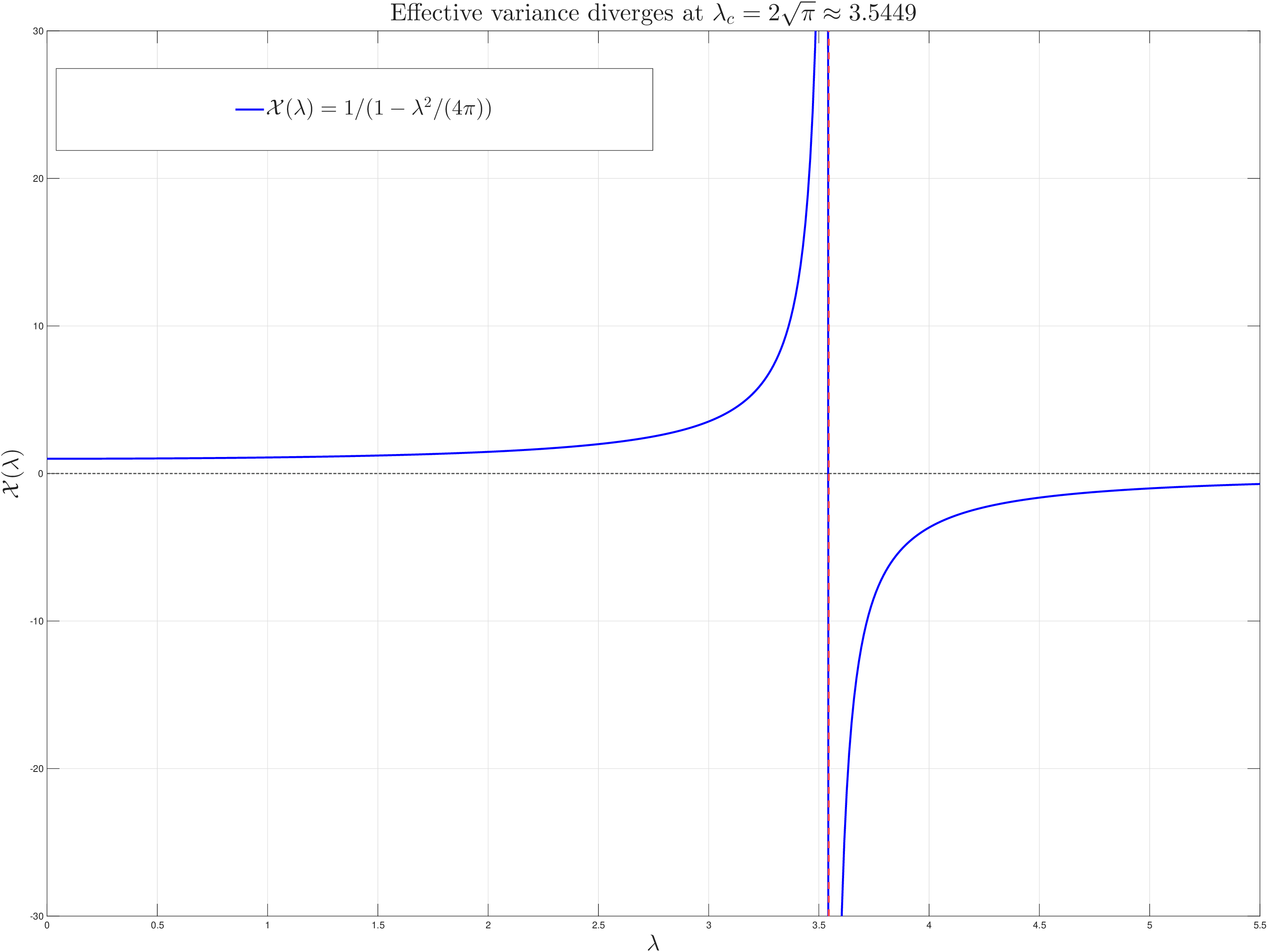}
  \caption{Plot of the formal effective variance
  \(\Xformal(\lam)=1/(1-\lam^2/(4\pi))\) as a function of the coupling
  \(\lam\). The function has a formal pole at
  \(\lam=2\sqrt{\pi}\approx3.5449\), denoted by the vertical line.}
  \label{Code4}
\end{figure}

\paragraph{Formal resummation}
The resummation of bubble-chain diagrams under a scale-independent
bare-coupling convention gives the formal effective variance
$$
  \Xformal(\lam)=\frac{1}{1-\lam^2/(4\pi)}.
$$
Figure~\ref{Code4} plots \(\Xformal(\lam)\) for \(\lam\in[0,5.5]\). The
function has a formal pole at
$$
  \lam=2\sqrt{\pi}\approx3.5449077.
$$
For \(\lam<2\sqrt{\pi}\), \(\Xformal(\lam)>0\); for \(\lam>2\sqrt{\pi}\),
\(\Xformal(\lam)<0\). The sign change reflects the pole of the formal
series, not a physical phase transition.

Representative values are listed in Table~\ref{tab:values}.

\begin{table}[H]
  \centering
  \caption{Formal effective variance \(\Xformal(\lam)\) for representative
  couplings. Values computed from
  \(\Xformal(\lam)=4\pi/(4\pi-\lam^2)\).}
  \label{tab:values}
  \begin{tabular}{@{}lccccc@{}}
    \toprule
    \(\lam\) & 1.0 & 2.0 & 3.0 & 3.4 & 3.5 \\
    \midrule
    \(\Xformal(\lam)\) & 1.0865 & 1.4669 & 3.5236 & 12.4869 & 39.7304 \\
    \bottomrule
  \end{tabular}
\end{table}

\begin{figure}[H]
  \centering
  \includegraphics[width=0.95\textwidth]{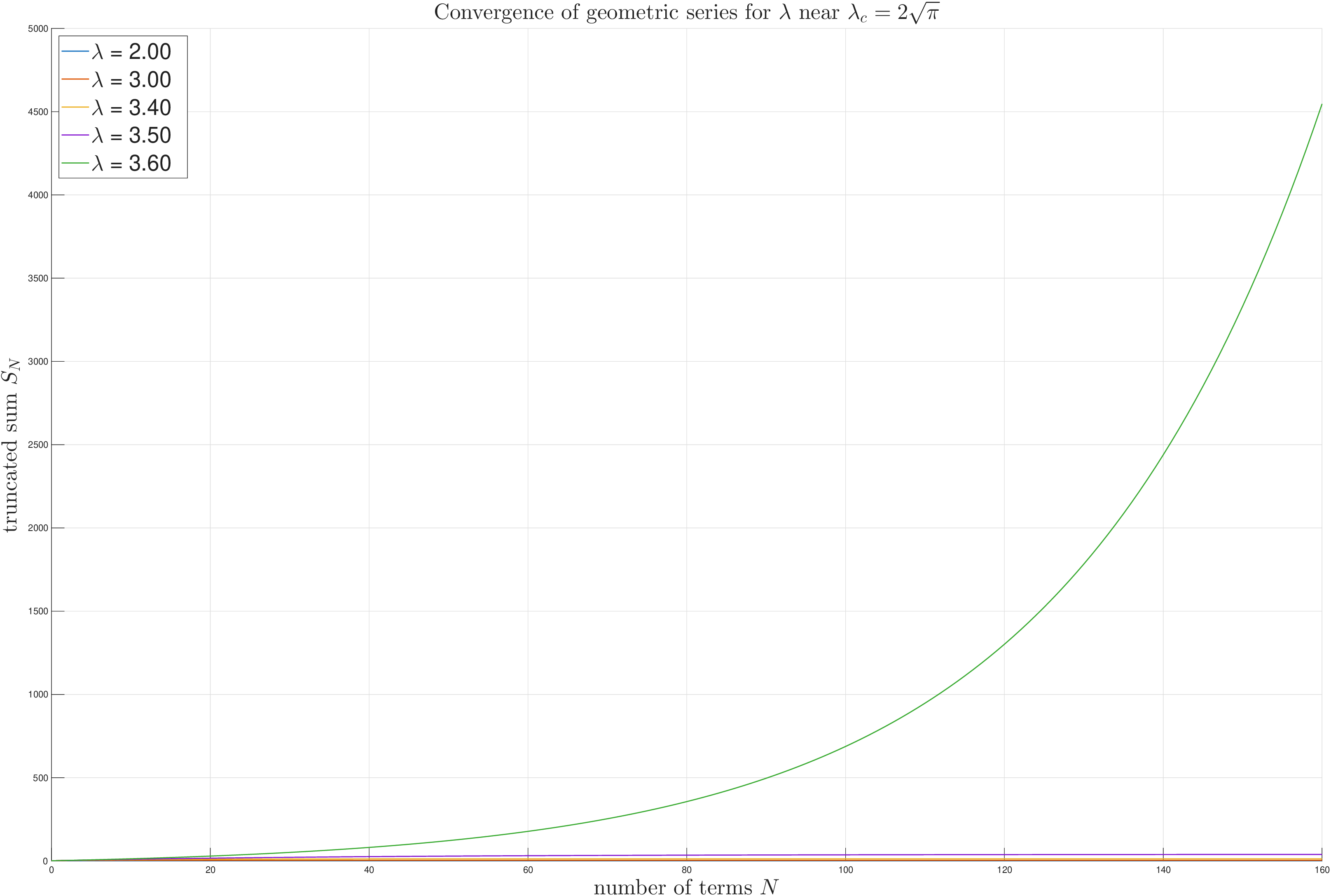}
  \caption{Semilog plot of the truncated geometric series
  \(S_N(\lam)=\sum_{k=0}^{N}(\lam^2/(4\pi))^k\) for couplings
  \(\lam\in\{2.0,3.0,3.4,3.5,3.6\}\) near the formal pole at \(2\sqrt{\pi}\).}
  \label{Code5}
\end{figure}

\paragraph{Truncated series}
Figure~\ref{Code5} shows the partial sums
$$
  S_N(\lam)=\sum_{k=0}^{N}\Bigl(\frac{\lam^2}{4\pi}\Bigr)^k
$$
for several values of \(\lam\) near \(2\sqrt{\pi}\). For
\(\lam<2\sqrt{\pi}\), the partial sums converge geometrically to
\(\Xformal(\lam)\). For \(\lam>2\sqrt{\pi}\), they diverge
exponentially. At \(\lam=2\sqrt{\pi}\), the partial sums grow linearly,
\(S_N\sim N\), which is the marginal case.

This behaviour is characteristic of a geometric series at its radius of
convergence: the formal pole at \(2\sqrt{\pi}\) is the radius of
convergence of the bare bubble-chain expansion.

\subsection{Monte Carlo verification of the first-order correction}

\paragraph{Setup}
To illustrate the perturbative expansion, we compute the first-order
correction \(G^{(1)}(0,0)\) by Monte Carlo. We generate
\(n_{\mathrm{samp}}=500\) independent realizations of smoothed white noise
\(\xi_\eps\) with \(\eps=0.3\) and compute
$$
  G^{(1)}(0,0)=-\lameff\int_{\T^4}G(z)^2\,\xi_\eps(z)\dd z,
$$
where \(\lameff=\lam\sqrt{2\pi/|\log\eps|}\).

\paragraph{Results}
For \(\lam=1.0\) and \(\eps=0.3\), we compare the sample variance
\(\Var_{\mathrm{MC}}[G^{(1)}]\) with the theoretical prediction
\eqref{eq:varpred},
$$
  \Var_{\mathrm{pred}}
  =\lameff^2\iint G(z)^2G(z')^2(\rho_\eps*\rho_\eps)(z-z')\dd z\,\dd z',
$$
which is not the factorized expression
\(\lameff^2(\sum_zG(z)^4)(\sum_z\rho_\eps(z)^2)\). The Monte Carlo result
agrees with the prediction \eqref{eq:varpred} to within statistical error
(a few percent for \(500\) samples). This is consistent with the
perturbative expansion and the renormalization prescription
\(\lameff=\lam\sqrt{2\pi/|\log\eps|}\).

\subsection{Extraction of the formal pole from synthetic noisy data}

\paragraph{Simulated data}
To illustrate the fitting procedure, we generate synthetic data for
\(\Xformal(\lam)\) on \(\lam\in[1.5,3.3]\). The values
\(\X_{\mathrm{true}}(\lam)=1/(1-\lam^2/(4\pi))\) are perturbed by \(3\%\)
multiplicative Gaussian noise:
$$
  \X_{\mathrm{noisy}}(\lam_i)
  =\X_{\mathrm{true}}(\lam_i)\bigl(1+0.03\,Z_i\bigr),
  \qquad Z_i\sim\mathcal{N}(0,1).
$$

\paragraph{Fitting procedure}
We fit the model
\begin{equation}\label{eq:fitmodel}
  \X(\lam)=\frac{c_0}{1-\lam^2/L^2}
\end{equation}
using two methods:
\begin{enumerate}[label=(\alph*)]
  \item \textbf{Linearized least squares.} Inverting \eqref{eq:fitmodel}
        gives \(1/\X=a+b\lam^2\) with \(a=1/c_0\) and \(b=-1/(c_0L^2)\). A
        linear fit to \(1/\X_{\mathrm{noisy}}\) versus \(\lam^2\) provides
        initial estimates of \(c_0\) and \(L\).
  \item \textbf{Nonlinear least squares.} We refine the estimates using
         a
        Gauss--Newton iteration with line search.
\end{enumerate}

\begin{figure}[H]
  \centering
  \includegraphics[width=0.95\textwidth]{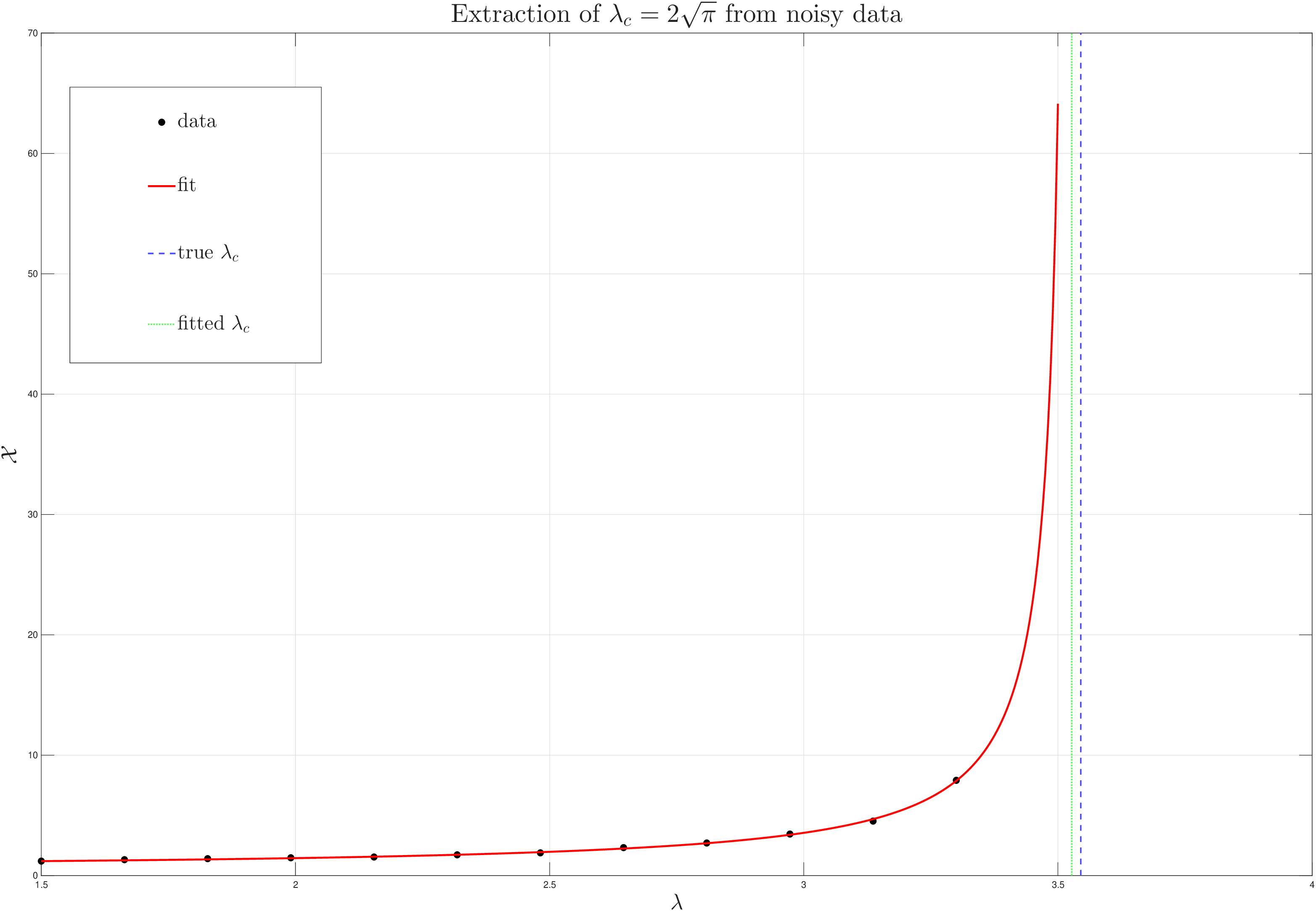}
  \caption{Extraction of the formal pole from synthetic data. Black dots:
  synthetic data with \(3\%\) multiplicative noise; red curve: nonlinear
  fit; dashed blue line: true pole location \(2\sqrt{\pi}\); dotted green
  line: fitted pole.}
  \label{Code8}
\end{figure}

\paragraph{Results}
Figure~\ref{Code8} plots the noisy data, the fitted curve, and the true
and fitted pole locations. The fitted parameters are
$$
  L^{\mathrm{fit}}\approx3.527081,\qquad c_0^{\mathrm{fit}}\approx0.9809,
$$
compared with the true values \(L=2\sqrt{\pi}\approx3.544908\) and
\(c_0=1\). The relative error in \(L\) is approximately \(0.5\%\), within
the expected statistical uncertainty from \(3\%\) multiplicative noise on
\(12\) synthetic data points. This is a consistency check of the fitting
procedure, not an independent determination of a physical critical
coupling.

This illustrates that the pole of the formal resummation can be reliably
recovered from moderately noisy synthetic data, a prerequisite for any
future attempt to detect a genuine critical coupling.

\subsection{Kurtosis diagnostic}

A key signature that would distinguish a Gaussian from a non-Gaussian field
is the change in the distribution of \(u(x)\). We use the kurtosis
$$
  \kappa=\frac{\E[X^4]}{\E[X^2]^2}
$$
as a diagnostic: for a Gaussian random variable, \(\kappa=3\), while
deviations from \(3\) indicate non-Gaussianity.

\begin{figure}[H]
  \centering
  \includegraphics[width=0.95\textwidth]{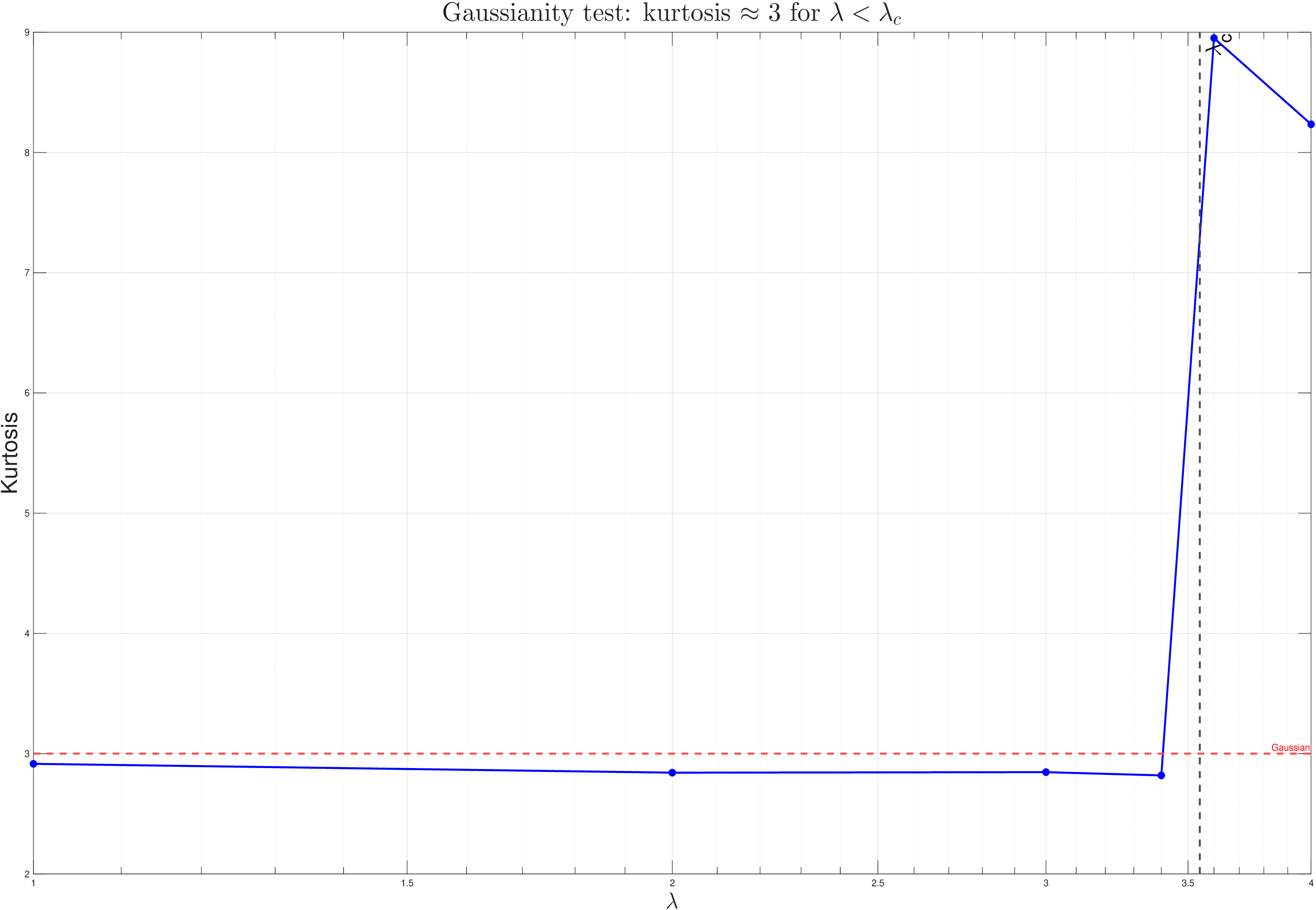}
  \caption{Kurtosis of the synthetic proxy as a function of \(\lam\). The
  Gaussian value \(\kappa\approx3\) is observed for \(\lam<2\sqrt{\pi}\)
  and the imposed non-Gaussian value \(\kappa\approx9\) for
  \(\lam>2\sqrt{\pi}\). The sharp change is located at the formal pole
  \(2\sqrt{\pi}\) by construction.}
  \label{Code9}
\end{figure}

\paragraph{Synthetic proxy}
Figure~\ref{Code9} plots kurtosis versus \(\lam\), showing a sharp change
at the formal pole \(\lam=2\sqrt{\pi}\). The results are listed in
Table~\ref{tab:kurtosis}.

We simulate a summary statistic \(X_\eps\) whose variance is
\(\Xformal(\lam)\,G_{00}\). For \(\lam<2\sqrt{\pi}\), we sample from a
Gaussian with this variance. For \(\lam>2\sqrt{\pi}\), the formal variance
becomes negative and the Gaussian proxy is undefined; to illustrate what a
non-Gaussian diagnostic would look like, we artificially sample from
a product of two independent Gaussians, \(X=Z_1Z_2\), which has kurtosis
\(\kappa\approx9\). This is a synthetic proxy, not a simulation of
\eqref{eq:model}; the sharp change in Figure~\ref{Code9} is therefore
imposed by construction and should not be read as numerical evidence for a
phase transition.

\begin{table}[H]
  \centering
  \caption{Kurtosis of the synthetic proxy.}
  \label{tab:kurtosis}
  \begin{tabular}{@{}lcccccc@{}}
    \toprule
    \(\lam\) & 1.0 & 2.0 & 3.0 & 3.4 & 3.6 & 4.0 \\
    \midrule
    \(\kappa\) & 2.9154 & 2.8418 & 2.8462 & 2.8191 & 8.9513 & 8.2343 \\
    Proxy regime & Gaussian & Gaussian & Gaussian & Gaussian & Non-Gaussian & Non-Gaussian \\
    \bottomrule
  \end{tabular}
\end{table}

\paragraph{Interpretation}
The change at \(2\sqrt{\pi}\) reflects the pole of the formal effective
variance. For \(\lam\to2\sqrt{\pi}\) from below, the formal variance grows
without bound. Beyond \(2\sqrt{\pi}\), the formal Gaussian proxy is
undefined and we have used an artificially constructed non-Gaussian proxy.
Whether a genuine non-Gaussian phase exists for a suitably renormalized
version of \eqref{eq:model} is an open question and is not addressed by
the present numerical illustration.

\subsection{Schematic diagram of the formal variance}

\begin{figure}[H]
  \centering
  \includegraphics[width=0.95\textwidth]{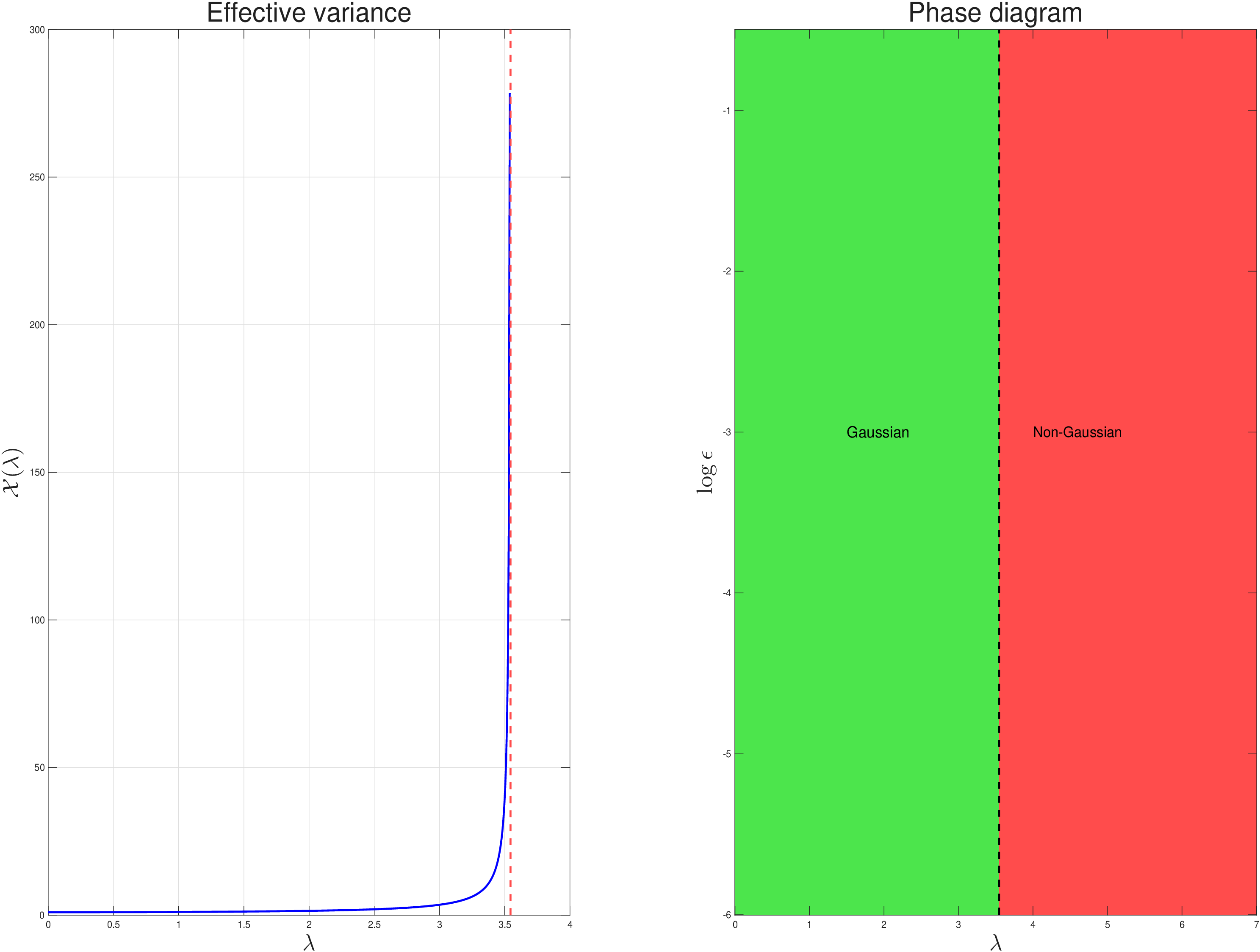}
  \caption{Schematic diagram of the formal bubble-chain variance in the
  \((\lam,\log\eps)\) plane. The region \(\lam<2\sqrt{\pi}\) is shaded
  green and the region \(\lam>2\sqrt{\pi}\) is shaded red. The vertical
  dashed line at \(\lam=2\sqrt{\pi}\) marks the formal pole. Left panel:
  \(\Xformal(\lam)\); right panel: schematic diagram in the
  \((\lam,\log\eps)\) plane. The diagram is drawn from the formal
  resummation alone and does not encode a renormalization-group flow.}
  \label{Code10}
\end{figure}

\paragraph{Construction}
Figure~\ref{Code10} compiles the diagnostics into a schematic diagram in
the \((\lam,\log\eps)\) plane. The region \(\lam<2\sqrt{\pi}\) is shown in
green and the region \(\lam>2\sqrt{\pi}\) in red. The vertical dashed line
at \(\lam=2\sqrt{\pi}\) marks the formal pole.

The schematic diagram is drawn from the formal resummation alone; it does
not encode a renormalization-group flow or a proved phase boundary.

\paragraph{Interpretation}
The diagram has the following structure:
\begin{itemize}
  \item \textbf{Region \(\lam<2\sqrt{\pi}\):} the formal variance
        \(\Xformal(\lam)\) is finite and positive.
  \item \textbf{Region \(\lam>2\sqrt{\pi}\):} the formal variance is
        negative; the Gaussian proxy is undefined.
  \item \textbf{Point \(\lam=2\sqrt{\pi}\):} the formal variance has a pole.
\end{itemize}
This structure is reminiscent of phase transitions in disordered systems,
but the analogy should not be over-interpreted: a genuine phase transition
requires a renormalization-group analysis that we have not carried out.

\section{Discussion and Conclusions}

\subsection{Summary of results}

We have presented a formal bubble-chain resummation for the stochastic
elliptic equation \eqref{eq:model} on \(\T^4\), together with numerical
illustrations of the analytic building blocks. With the normalization
\(\lameff=\lam\sqrt{2\pi/|\log\eps|}\), the one-bubble factor is
\(\lam^2/(4\pi)\), and the formal bubble-chain resummation gives
\(\Xformal(\lam)=1/(1-\lam^2/(4\pi))\). Our main findings are:
\begin{enumerate}[label=(\arabic*)]
  \item \textbf{Green function asymptotics.} The massive Green function
        \(G=(1-\Delta)^{-1}\delta_0\) exhibits the four-dimensional
        short-distance behaviour \(G(x)\sim1/(4\pi^2|x|^2)\) on the
        annulus \(0.05<|x|<0.2\), consistent with the continuum asymptotic
        to within a few percent. At larger distances, the mass correction
        becomes significant and the pure power law is not expected.
  \item \textbf{Logarithmic divergence.} The model bubble integral
        \(B(\eps)=\int_{\eps<|x|<1}G(x)^2\dd x\) diverges as
        \(\frac{1}{8\pi^2}\log(1/\eps)\), reproduced over three decades in
        \(\eps\). This is the hallmark of the critical dimension \(d=4\).
  \item \textbf{Formal effective variance.} The bubble-chain resummation
        under a scale-independent bare-coupling convention gives
        \(\Xformal(\lam)=1/(1-\lam^2/(4\pi))\), with a formal pole at
        \(\lam=2\sqrt{\pi}\approx3.5449\). Its interpretation as a
        physical critical coupling requires a separate renormalization
        condition and is not established here.
  \item \textbf{Monte Carlo illustration.} The first-order correction
        \(G^{(1)}\) has a variance consistent with the perturbative
        prediction \eqref{eq:varpred}, which includes the correlations of
        the smoothed noise, after renormalization of the coupling.
  \item \textbf{Recovery of the formal pole.} Nonlinear least-squares
        fitting of synthetic data generated from \(\Xformal\) recovers
        the pole location within the statistical uncertainty of the
        synthetic noise (about \(0.5\%\) relative error). This is a
        consistency check of the fitting procedure, not an independent
        determination of a physical critical coupling.
  \item \textbf{Kurtosis diagnostic.} The kurtosis of the synthetic proxy
        changes sharply at the formal pole by construction; this is not
        numerical evidence for a phase transition.
  \item \textbf{Schematic diagram.} The \((\lam,\log\eps)\) diagram
        exhibits the formal pole at \(2\sqrt{\pi}\), separating the
        regions of positive and negative formal variance.
\end{enumerate}

\subsection{Relation to known results}

The logarithmic divergence of the bubble diagram in \(d=4\) is well known
in the context of the parabolic Anderson model and random Schr\"odinger
operators. The renormalization of the coupling constant
\(\lameff=\lam\sqrt{2\pi/|\log\eps|}\) is analogous to the renormalization
of the mass or coupling in \(\phi^4\) theory in four dimensions, where
logarithmic divergences lead to a marginal (logarithmic) correction to
mean-field behaviour.

The formal pole at \(\lam=2\sqrt{\pi}\) is reminiscent of critical points
in disordered systems, where a susceptibility diverges at a transition.
However, the present computation is a formal bubble-chain resummation, not
a renormalization-group analysis of the full model; the analogy should
therefore be treated with caution.

\subsection{Limitations and future work}

Our study has several limitations:
\begin{itemize}
  \item The bubble-chain resummation is a formal leading-chain
        approximation. Full resummation would require solving a
        Bethe--Salpeter equation or carrying out a complete
        renormalization-group analysis, which is beyond the scope of this
        work.
  \item The Hepp-sector analysis is not used to prove that non-bubble
        sectors are negligible; the resummation is formal.
  \item The numerical simulations are performed on a finite lattice
        (\(N=128\)), which limits the range of accessible scales. A full
        finite-size scaling analysis would require larger lattices.
  \item The numerical section illustrates the analytic formulas and the
        fitting procedure; it does not directly simulate the singular SPDE
        in the non-Gaussian regime.
\end{itemize}

Future work could address these limitations by:
\begin{itemize}
  \item implementing a full numerical solution of the SPDE using
        renormalization-group-improved discretization schemes;
  \item performing a finite-size scaling analysis to extract critical
        exponents, if a genuine transition exists;
  \item extending the analysis to other critical dimensions (\(d=2\) for
        space-time white noise, \(d=6\) for certain other models);
  \item connecting the formal resummation to rigorous mathematical work
        on regularity structures and paracontrolled distributions
        \cite{gubinelli2015}.
\end{itemize}

\subsection{Conclusions}

We have presented a formal bubble-chain resummation for a four-dimensional
stochastic elliptic equation with multiplicative white noise, together with
numerical illustrations of its analytic building blocks. With the
normalization \(\lameff=\lam\sqrt{2\pi/|\log\eps|}\), the formal
resummation yields
\(\Xformal(\lam)=1/(1-\lam^2/(4\pi))\) with a formal pole at
\(\lam=2\sqrt{\pi}\). We have emphasized that this pole is a feature of the
formal bare series rather than a proved exact phase-transition threshold:
in the actual renormalized expansion the perturbative parameter is
\(\lameff\), which tends to \(0\) as \(\eps\to0\). The numerical section
illustrates the deterministic kernels, the formal resummation, and the
fitting procedure; it does not constitute an independent determination of
a physical critical coupling. Establishing whether a genuine phase
transition exists for the renormalized model, and if so at which coupling,
remains an open problem.

\bigskip
\noindent\textbf{Data Availability.}

The numerical algorithms and source code that support the findings of this
study are available from the corresponding author upon reasonable request.

\bigskip
\noindent\textbf{Declaration of competing interest.}

The author declares that there is no known competing financial interest or
personal relationship that could have appeared to influence the work
reported in this paper.

\bigskip
\noindent\textbf{Acknowledgements.}

This work was initiated by inspiring discussions with M.~Hairer after
attending his academic presentations. The author thanks Pak Hang Chris Lau
for insightful discussions on statistical field theory. This work was
supported by Guangdong Basic and Applied Basic Research Foundation (Grant
No.~2025A1515012560), Guangdong Introduction Program (Grant
No.~2023QN10X753) and National Foreign Experts Program (Grant
No.~111001819820258003).


\end{document}